\documentclass[a4paper]{amsart}

\usepackage[all]{xy}

\usepackage{amsmath,amssymb}
\usepackage{amscd,eucal,amsthm, eufrak}
\usepackage{epsfig}
\CompileMatrices

\newtheorem{theorem}{Theorem}[section]
\newtheorem{lemma}[theorem]{Lemma}
\newtheorem{proposition}[theorem]{Proposition}
\newtheorem{corollary}[theorem]{Corollary}
\theoremstyle{definition}
\newtheorem{definition}[theorem]{Definition}
\newtheorem{example}[theorem]{Example}

\newtheorem{remark}[theorem]{Remark}

\numberwithin{equation}{theorem}
\def\id{{\rm id}}

\def\Ker{{\rm Ker}\,}

\newcounter{itemnumber}

\begin{document}

\sloppy

\title[ Variation of Mumford's  quotients of  a flag variety
 ] { Variation of Mumford's  quotients for
the maximal torus action on a flag variety }

\author[V.~Zhgoon ]{Vladimir S. Zhgoon }
\thanks{}
\address{Department of Higher Algebra,
Faculty of Mechanics and Mathematics, Moscow State Lomonosov
University, Vorobievy Gory, GSP-2, Moscow, 119992, Russia}
\email{zhgoon@mail.ru}

\begin{abstract}
We study a variation of  Mumford's quotient for the action of a
maximal torus $T$ on a flag variety $G/B$ depending on projective
embedding
 $G/B \hookrightarrow \Bbb P(V(\chi))$,
where  $T$-linearization is induced by the standard
$G$-linearization. We describe the linear spans of the supports of
 semistable orbits, that allow us to calculate the rank of the
Picard group of  quotient $(G/B)^{ss}/\!\!/T$, when $G$ does not
contain simple factors of  type $A_n$.

\end{abstract}

\maketitle
 Let  $G$  be a semisimple algebraic group over an
algebraically closed field of characteristic zero, $T$ a maximal
torus in $G$, and $B$ a Borel group containing $T$.
 Consider the action of $T$ on $G/B$ by left multiplication.
Let $\chi$  be a strictly dominant weight. It is well known that
$G/B$ can be embedded  $G$-equivariantly in the projectivization
$\Bbb P(V(\chi))$ of à simple module $V(\chi)$ of the highest
weight $\chi$ as the projectivization of the orbit of the highest
weight vector. All $G$-equivariant embeddings $G/B\hookrightarrow
\Bbb P^N$ can be obtained by this construction. Denote by
$L_{\chi}$ the restriction on $G/B$ of the $G$-linearized sheaf
${\mathcal O} (1)$ on $\Bbb P(V(\chi))$. According to ~\cite
{mum}, a Zariski open subset $X^{ss}_{L_\chi}$ of the flag variety
$X=G/B$ can be defined, in such way that there exists a
categorical quotient $X^{ss}_{L_\chi}/\!\!/T$ for the torus
action. In this paper we study  variation of the Mumford's
quotient depending on the $T$-linearized sheaf $L_{\chi}$. We also
describe the linear spans of the supports of  semistable
$T$-orbits, that allow us to calculate the rank of the Picard
group $Pic(X^{ss}_{L_\chi}/\!\!/T)$ of the projective variety
$X^{ss}_{L_\chi}/\!\!/T$ (depending on the strictly dominant
weight $\chi$). We note that in this case
$Pic(X^{ss}_{L_\chi}/\!\!/T)$ is finitely generated free abelian
group.

For the convenience  of the reader we remind the definition of the
set of (semi)stable points.

\begin{definition} Let $X$  be an algebraic variety with an
action of $G$, and $L$
 be an invertible ample $G$-linearized sheaf on $X$.

(i) Following Mumford we define the set  of semistable points as
$$X_L^{ss} =\{ x\in X :\exists n >0, \exists \sigma \in
\Gamma(X,L^{\otimes n})^G,
 \sigma(x)\not = 0\}.$$
(ii) The set of  stable points is defined as
 $$ X_L^s=\{ x \in X_L^{ss} :  \text {the orbit} \ Gx \ \text{is closed}
\  X_L^{ss}  \ \text {and the stabilizer}  \  G_x  \ \text {is
finite} \}.$$
\end{definition}

The orbits of maximal torus  $T$ on the flag varieties were
studied for instance in the papers ~\cite {Carrel Kurth},~\cite
{dabr},~\cite {sen}. In ~\cite {dabr} the normality of the
closures of typical $T$-orbits on $G/P$ was proved by R.Dabrowski
(where $P\supseteq B$
--- is a parabolic subgroup). In ~\cite {Carrel Kurth}
the normality of the closures of non-typical $T$-orbits on $G/P$
was studied by J.B.Carrell and A.Kurth. In ~\cite {sen}
S.Senthamarai Kannan has found all the flag varieties $G/P$ with
the property that the equality $(G/P)^{ss}_L=(G/P)^{s}_L$ is
satisfied for an invertible sheaf $L$.

{\vspace {1ex}}

In the first part we derive the Seshadri criterion {\cite
{seshadri}} of the semistability of the point on $G/B$ for the
maximal torus action and decompose the Weyl chamber $C$ in the
$\rm GIT$-equivalence classes, for the points with the same sets
$X^{ss}_{L_\chi}$. In the second part we describe the linear spans
of the supports of the   $T$-orbits of the subvariety
$X^{ss}_{L_\chi}$. In the third part we apply the previous results
to the calculation of the rank of $Pic(X^{ss}_{L_\chi}/\!\!/T)$
when the group $G$ does not contain simple factors of type $A_n$.

The author is grateful to his scientific advisor  I.V.Azhantsev
 for posing the problem, and for the constant attention to the work. He also wants
to thank  D.A.Timashev, whose remarks leaded to the simplification
of the proof of  Theorem 3.1, and E.B.Vinberg for the valuable
comments and discussing several questions.

In the subsequent work the author will consider the case  $A_n$.
There is  also a hope to study the variation of Mumford's quotient
not only for standard, but for all possible $T$-linearizations of
 ample line bundles on $G/B$.

{\vspace {1ex}}

{\bf Notations and conventions.}

{\vspace {1ex}}

 By gothic letters we  denote Lie algebras corresponding to  Lie groups.

 $\Xi=\Xi(T)$ --- the lattice of  characters of $T$. Its dual is
identified with the lattice of one parameter subgroups

  $\Lambda(T)$ by the pairing invariant under the action of Weyl group, that we denote by $\langle .; .\rangle$.

 $\Xi_{\Bbb Q}=\Xi \otimes _{\Bbb Z} \Bbb Q$ ---  rational characters of the torus $T$.

 $\Lambda_{\Bbb Q}=\Lambda \otimes _{\Bbb Z} \Bbb Q$ ---  rational one parameter subgroups of the torus $T$.

 $W=N_G(T)/T$ --- Weyl group.

 On $\Xi_{\Bbb Q}$ we have  $W$-invariant scalar product
 $(.;.)$, that we use to identify $\Xi_{\Bbb Q}$
with $\Lambda_{\Bbb Q}$.

 $\Delta  $ is the root system of the Lie algebra $\mathfrak g$ corresponding to $T$.
 $\Delta ^{+} (\Delta ^{-})  $ is the system of the positive (negative) roots corresponding to the Borel subalgebra  $\mathfrak b \subset \mathfrak g$.
 $\Pi$ is the system of the simple roots.
$w_0\in W$ is the longest element in the Weyl group.
 $C$ is the positive Weyl chamber.

 Let $\widetilde{\Delta}\subset \Delta$ be the subset of roots that is an abstract system of roots.
 $C_{\widetilde{\Delta}}:=\{ \chi \in \Xi_{\Bbb
Q}|(\chi;\alpha_i)\geqslant 0$ for all $\alpha_i \in
\widetilde{\Delta}^+\}
 $ denotes the positive Weyl chamber for the root system
$\widetilde{\Delta}$.

We call the root subsystem $\widetilde{\Delta}\subset \Delta$
saturated, if the following property holds:
$\widetilde{\Delta}=\langle\widetilde{\Delta}\rangle\cap\Delta$ (
here $\langle\widetilde{\Delta}\rangle$ denotes the linear span of
the root system $\widetilde{\Delta}$).

 Let $V(\chi)$ be a simple $G$-module with the highest strictly dominant
weight $\chi\in C\cap \Xi$ and the highest weight vector $v_{\chi}
\in V(\chi)$.

 Consider the action of $T$ on the linear space $V=\bigoplus \limits_{\lambda \in \Xi} V_{\lambda}$
($V_\lambda$
 is the weight component of the weight  $\lambda \in \Xi$), $v=\sum \limits_{\lambda \in \Xi} v_{\lambda}$,
where
  $v_{\lambda}\in V_\lambda $ and $v_{\lambda}\neq 0$. Denote by $supp(v) \subset
 \Xi_{\Bbb Q}$ the convex hull of the weights of vector $v$. One may
notice that
 $supp(v)=supp(tv)$ for every $t \in T$, so the support of a $T$-orbit is correctly defined, that we  denote it by $supp(Tv)$.
We  also define the support of the $T$-invariant subset as the
convex hull of the supports of all $T$-orbits from this set. Let
$x=\langle v\rangle\in \Bbb P(V)$ be a point corresponding to  the
vector $v\in V$. Denote by $supp(Tx \hookrightarrow \Bbb P(V))$
the support of the orbit $Tv \hookrightarrow V$. The set of
weights of the  $T$-orbit in the embedding $Tx\subset\Bbb
P(V(\chi))$ is denoted by $Pd_{\chi}(Tx)$.

 If it is not stated otherwise, everywhere in the article we consider the flag
 variety $G/B$ with the
 $G$-equivariant embedding in $ \Bbb P(V(\chi))$
 as a close $G$-orbit of the line, generated by the highest weight vector $v_{\chi}$.

{\vspace {1ex}}

 \section{A criterion of stability for a point  in $G/B$}

{\vspace {1ex}}

The stability of the point depends on the Shubert variety to which
it belongs. It is essential to study the geometry of the
projective embedding  for the Schubert variety.  Thus we have to
remind the lemma from the paper of
 I.N.Berstein, I.M.Gelfand, S.I.Gelfand, that  describes the structure
of the Schubert variety under given projective embedding.

{\vspace {2ex}}

\begin{lemma}{~\cite[2.12] {BGG}}  Let $w \in W$ be an
element of the Weyl group, $BwB/B$ the
 corresponding Schubert cell.
 Consider the closed embedding $G/B \hookrightarrow \Bbb
P(V(\chi))$, where $\chi$ is a strictly dominant weight. Let $f\in
V({\chi})$ be a vector from the orbit of the highest weight
vector. Then $\langle f\rangle \in BwB/B$ iff $w\chi \in supp(f)$
and $f \in \mathfrak U(\mathfrak b)v_{w\chi}$, where $\mathfrak
U(\mathfrak b)$ is the universal enveloping algebra of the Lie
algebra $\mathfrak b$, and $v_{w\chi}=w v_{\chi}$ is a vector of
the weight ${w\chi}$, entering the irreducible representation
$V(\chi)$ with  multiplicity one.
\end{lemma}

Applying this lemma it is easy to prove the following
semistability criterion due to Seshadri ~\cite[Prop.~1.5]
{seshadri} . We  give the proof that differs from the original
one.  The ideas similar to those of the proof will be used later.
 {\vspace {1ex}}
But  first we remind the definition of  Mumford's numerical
function for a torus action and   Mumford's numerical criterion of
stability.

\begin{definition} Let $L$ be an ample $T$-linearized line
bundle on a $T$-variety $X$ defining the $T$-equariant embedding
of $X$ in the projective space $\Bbb P(V)$. Let $\lambda \in
\Lambda(T)$ be a one-parameter subgroup. For a point $x\in \Bbb
P(V)$ we calculate  Mumford's numerical function by the formula
$$\mu^{L}(x,\lambda)=\min\limits_{\tau \in
supp(Tx)}\langle\tau;\lambda\rangle.$$

\end{definition}

\begin{proposition}{ Mumford's numerical criterion (~\cite
{mum})} Let $X$ be a variety with a $T$-action   and $L$ be an
ample $T$-linearized line bundle defining $T$-equivariant
embedding $X$ in the projective space $\Bbb P(V)$.  The point
$x\in \Bbb P(V)$ is (semi)stable iff
$\mu^{L}(x,\lambda)(\leqslant)<0$ for every nontrivial
one-parameter subgroup $\lambda \in \Lambda(T)$.
\end{proposition}
{\vspace {2ex}}

 \begin{proposition} {(Seshadri ~\cite {seshadri})} Let $C$ be the Weyl chamber, $x \in G/B$
and $x=bwB/B$. Consider a very ample line bundle $L_{\chi}$,
corresponding to the strictly dominant weight $\chi$. Let $\lambda
\in C$ be a one-parameter group belonging to the Weyl chamber.
Then we have $\mu^{L_{\chi}}(x,\lambda)=\langle w\chi;\lambda
\rangle $.
\end{proposition}

\begin{proof} Consider the support of the orbit $Tx$. As $Tx
\subset BwB/B $,  applying the previous lemma, we have $w\chi \in
supp(Tx)$. Also we have $supp(Tx) \subset
supp(BwB/B)=supp(\mathfrak U(\mathfrak b)v_{w\chi}) $, where the
last equality is again the consequence of the previous lemma.

If the weight $\tau$ belongs to $supp(\mathfrak U(\mathfrak
b)v_{w\chi})$, then $\tau =w\chi+\sum\limits_{\alpha_i \in
\Delta^+}c_i\alpha_i $, where $c_i \geqslant 0$. Considering the
pairing of $\tau$ with the one-parameter subgroup  $\lambda$ we
get:
$$ \langle\tau;\lambda\rangle=\langle w\chi;\lambda\rangle+\sum\limits_{\alpha_i
\in \Delta^+}c_i\langle\alpha_i;\lambda\rangle \geqslant \langle
w\chi;\lambda\rangle,$$

\noindent as  $c_i \geqslant 0$ and
$\langle\lambda;\alpha_i\rangle \geqslant 0$. Then we may derive
the expression for the numerical function:
$$\mu^{L_{\chi}}(x,\lambda)=\min\limits_{\tau \in
supp(Tx)}\langle\tau;\lambda\rangle=\langle
w\chi;\lambda\rangle.$$
\end{proof}

{\vspace {1ex}} Now we are ready to get the description of the
semistable point set. We need to introduce the following
definition.

{\vspace {1ex}}

\begin{definition} We call $w \in W$ $\chi$-semistable
 if $\langle w\chi;\lambda\rangle\leqslant 0$ for every $\lambda \in C$.
 We define the set  of all such elements by $W^{st}_{\chi}$.
\end{definition}
{\vspace {1ex}}

\begin{theorem} Consider   $G$-equivariant closed embedding
$G/B \hookrightarrow \Bbb P(V(\chi))$. Then we can find the  set
of semistable points of the  action of  torus $Ò$(with the
linearization coming from the standard $T$-action  on $V(\chi)$)
via the following formula:
$$ X^{ss}_{L_{\chi}}=\bigcap \limits_{\tilde{w}\in W}\bigcup \limits_{w \in W^{st}_{\chi}}\widetilde{w}BwB/B $$
\end{theorem}

\begin{proof} Indeed,  $\bigcup \limits_{w \in
W^{st}_{\chi}}BwB/B$
 is the set of such $x$, that
$\mu^{L_{\chi}}(x,\lambda)\leqslant 0 $ for every one-parameter
subgroup $\lambda \in C$.

By the well known equality for numerical functions one has
$\mu^{L_{\chi}}(x,\lambda)=\mu^{L_{\chi}}(\widetilde{w}x,\widetilde{w}\lambda)$,
where $\widetilde{w}$ is a representative of the Weyl group
element in the normalizer of $T$.

From the above we can rewrite the condition of the semistability
of $x$:  $\mu^{L_{\chi}}(x,\lambda)\leqslant 0$ for every $\lambda
\in \Lambda_{\Bbb Q}(T)$ as
$\mu^{L_{\chi}}(\widetilde{w}x,\lambda)\leqslant 0$ for every
$\lambda \in C$ and $\widetilde{w} \in W$.  By this we get  the
following formula for the set satisfying these conditions:$$
X^{ss}_{L_{\chi}}=\bigcap \limits_{\tilde{w}\in
W}\widetilde{w}\Bigl(\bigcup \limits_{w \in W^{st}_{\chi}}BwB/B
\Bigr) =\bigcap \limits_{\tilde{w}\in W}\bigcup \limits_{w \in
W^{st}_{\chi}}\widetilde{w}BwB/B. $$

\end{proof}

\begin{remark} The set $X^{ss}_{L_{\chi}}$ is the largest
$N_G(T)$-invariant subset in $\bigcup \limits_{w \in
W^{st}_{\chi}}BwB/B $.
\end{remark}

\begin{proposition} Consider two ample $T$-linearized line
bundles $L_{\chi_1}$ and $L_{\chi_2}$. Then
$X^{ss}_{L_{\chi_1}}=X^{ss}_{L_{\chi_2}}$ implies
$W^{st}_{\chi_1}=W^{st}_{\chi_2}$.
\end{proposition}

\begin{proof} It is sufficient to prove that for every element
$w\in W^{st}_{\chi}$ we may find in $BwB/B$ a semistable
$T$-orbit. It will imply that the set $W^{st}_{\chi}$ is uniquely
defined by $X^{ss}_{L_{\chi}}$. We have $\langle
w\chi;\lambda\rangle\leqslant 0$ for every $\lambda \in C$, that
is  equivalent to $0 \in w\chi+\sum\limits_{\alpha_i \in
\Delta^+}\Bbb Q_+\alpha_i$. So by  Corollary 2.6 (that we shall
prove in the next part because of its technical difficulty) $0\in
supp(BwB/B)$. For typical orbit  $Tx$ from the Schubert cell
$BwB/B$ we have $supp(Tx)=supp(BwB/B)$, so $0\in supp(Tx)$, that
means that a typical orbit from the cell $BwB/B$ is semistable.
\end{proof}

We shall decompose the Weyl chamber in cells $C_i$
---  $\rm GIT$-equivalence classes characterized by the following property: two characters $\chi_1,\chi_2$
belong to the same cell $C_i$ iff
$X^{ss}_{L_{\chi_1}}=X^{ss}_{L_{\chi_2}}$, that is equivalent to
$W^{st}_{\chi_1}=W^{st}_{\chi_2}$ by the preceding proposition.

{\vspace {1ex}}

Let $A$ be the cone generated by the simple roots. For an element
$\chi \in C^0$ in the relative interior of the Weyl chamber we
define the cone $\sigma_{\chi}=C\cap\bigcap \limits_{\stackrel
{w\in W} {\chi \in wA}} wA$.

 \begin{theorem} (Variation of  Mumford's quotients)

The set of the cones $\{\sigma_{\chi}| \chi \in C^0\}$ is finite.
They form the fan with the support $C$. The internal points of
these cones correspond to the $\rm GIT$-equivalence classes.
\end{theorem}
\begin{proof} The finiteness of the set of  cones follows from
the fact, that we get $\sigma_{\chi}$ as  intersections of the
cones from the finite set $\{w_iA | w_i\in W\}$. It is easily seen
from the definition, that for the cones $\{\sigma_{\chi}| \chi \in
C^0\}$ we have two cases:  the face
 of one cone doesn't contain  internal points of the other, one
 cone is the face of the other ( in particular these cones may coincide).

According to the definition, $w\in W^{st}_{\chi}$ when $\langle
w\chi;\lambda\rangle\leqslant 0$ for every $\lambda \in C$. This
is the same as
 $w\chi \in -A$, which is equivalent to $\chi \in -w^{-1} A$.
So the internal points of the cone  $\sigma_{\chi}$ are in one to
one correspondence with the elements  $\widetilde{\chi}$, for
which $W^{st}_{\widetilde{\chi}}=W^{st}_{\chi}$. \end{proof}
{\vspace {2ex}}

 Let us begin the study of the set of semistable points . In the subsequent
propositions we  prove that for semisimple groups not containing
 simple components of  type  $A_n$, the codimension of the
set of non-stable points is strictly greater than one.

{\vspace {2ex}}

\begin{lemma} Assume that the group  $G$ does not contain
the simple components of type $A_n$. Let $s_{\alpha} \in W$
 be the reflection corresponding to a root $\alpha$. Then $s_{\alpha}w_0$
belongs to the set $W^{st}_{\chi}$, for any $\chi \in C^0$.
\end{lemma}

\begin{proof} We have to check that $\langle s_{\alpha}w_0
\chi;\lambda\rangle\leqslant 0$ for every $\lambda \in C$. Denote
by $\pi_\alpha$ the fundamental weight dual to the simple root
$\alpha $. If $\beta,\gamma \in \Pi $, then
$(\pi_{\gamma};\beta)=\frac{(\beta;\beta)}{2}\delta_{\gamma
\beta}$. Also we have the equality $s_{\gamma}
\pi_{\gamma}=\pi_{\gamma}-\gamma$.

{\vspace {1ex}}

 The weight $\chi$ is dominant, so  $w_0 \chi$ is antidominant $w_0 \chi = -\sum  c_{\beta}
\pi_{\beta}$, where $c_\beta > 0$. As $\lambda \in C$ we have
$\lambda =\sum a_\beta \pi_\beta$, where $a_\beta \geqslant 0$. So

$-\langle s_{\alpha}w_0 \chi;\lambda\rangle=\langle s_{\alpha}\sum
c_\beta \pi_\beta;\sum a_\beta
\pi_\beta\rangle=(c_{\alpha}(\pi_{\alpha}-\alpha)+ \sum
\limits_{\beta \neq \alpha ;\beta \in
\Pi}c_{\beta}\pi_{\beta};a_{\alpha}\pi_{\alpha}+\sum
\limits_{\beta \neq \alpha; \beta \in \Pi}a_{\beta}\pi_{\beta})=
c_{\alpha}a_{\alpha}((\pi_{\alpha};\pi_{\alpha})-\frac{(\alpha;\alpha)}{2})+(\sum
\limits_{\beta \neq \alpha ;\beta \in
\Pi}c_{\beta}\pi_{\beta};\sum \limits_{\beta \neq \alpha; \beta
\in \Pi}a_{\beta}\pi_{\beta})+(c_{\alpha}\pi_{\alpha};\sum
\limits_{\beta \neq \alpha; \beta \in
\Pi}a_{\beta}\pi_{\beta})+(\sum \limits_{\beta \neq \alpha ;\beta
\in \Pi}c_{\beta}\pi_{\beta};a_{\alpha}\pi_{\alpha})\geqslant 0$.

The last equality holds since
$(\pi_{\alpha};\pi_{\alpha})-\frac{(\alpha;\alpha)}{2}\geqslant 0$
is true for all simple root systems distinct from $A_n$ (that is
not difficult to get from Table 2  ~\cite{Vin}), other summands
are not negative since they are scalar products of
 two dominant weights.
\end{proof}

{\vspace {1ex}}

\begin{theorem} Let us assume that the group  $G$ does not
contain the simple components of  type $A_n$. Then the complement
to the semistable point set in $G/B$ is of codimension strictly
greater than one.
\end{theorem}

\begin{proof} We shall use the formula describing the set of
semistable points  taking into account that, by the previous
lemma, $s_{\alpha}w_0 \in W^{st}_{\chi}$ for every $\alpha \in
\Pi$ and
 $w_0 \in W^{st}_{\chi}$:

$$
X^{ss}_{L_{\chi}}=\bigcap \limits_{\tilde{w}\in
W}\widetilde{w}\bigcup \limits_{w \in W^{st}_{\chi}}BwB/B \supset
\bigcap \limits_{\tilde{w}\in W}\widetilde{w}( \bigcup
\limits_{\alpha \in \Pi}Bs_{\alpha}w_0B/B \cup Bw_0B/B)
$$

 The set $\bigcup \limits_{\alpha \in \Pi}Bs_{\alpha}w_0B/B
\cup Bw_0B/B$ has already the complement of codimention 2 that
implies that  $X^{ss}_{L_{\chi}}$ has the complement of
codimension not smaller than 2, as it is equal to the finite
intersection of sets with the complement of codimention not
smaller than 2 (translation by $\widetilde{w}$ does not change the
codimension of the complement).
\end{proof}
{\vspace {1ex}}

\begin{corollary} Under the conditions of
Theorem 1.11, ${\rm Pic}(G/B)={\rm Pic}(X^{ss}_{L_{\chi}})$ and
${\rm Pic}_T(G/B)={\rm Pic}_T(X^{ss}_{L_{\chi}})$, where ${\rm
Pic}_T(X)$ denotes the group of $T$-linearized line bundles on the
$T$-variety $X$.
\end{corollary}

\begin{example} We give an example which shows that for
the root system  $A_2$ the statements of  Lemma 1.10 and Theorem
1.11 are not true.

Let $\chi =a \pi_1+ b\pi_2$, where $a>b$ (pic.1). Then $\langle
s_{\alpha_2}w_0\chi;\lambda\rangle >0$ for $\lambda=\pi_2\in C$.
So $s_{\alpha_2}w_0\notin W^{st}_{\chi}$ and the divisor
$Bs_{\alpha_2}w_0B/B$ consists of unstable points.
\end{example}
\begin{center}
  \epsfig{file=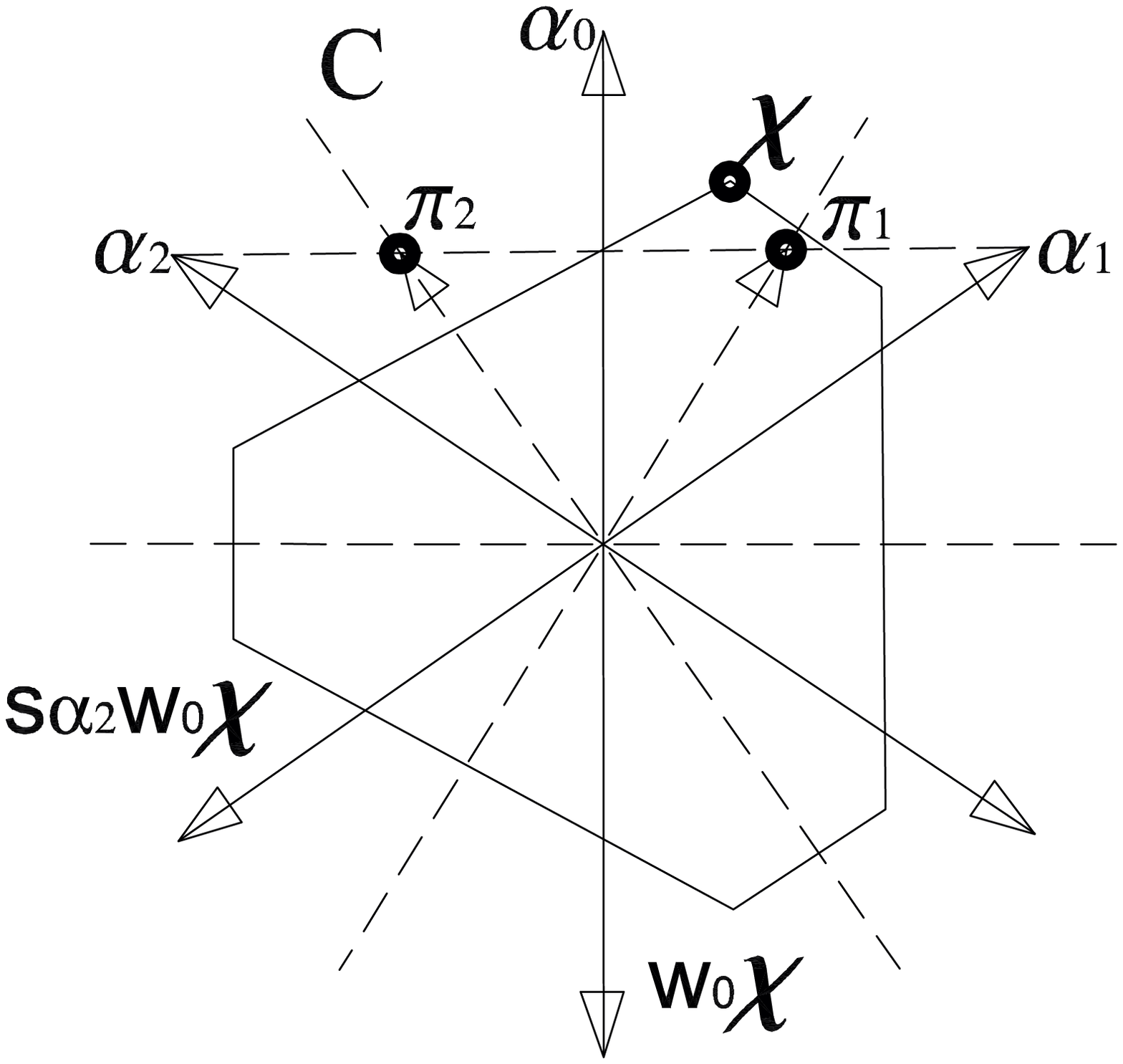,height=5cm,angle=90,clip=}
pic.1
\end{center}

\section{ The linear spans of  supports for  semistable
$T$-orbits}

In this part of the work we describe the linear spans of the
supports for semistable $T$-orbits on  $X^{ss}_{L_{\chi}}$.

{\vspace {1ex}}

 Let $x \in G/B$ be written in the form $x=bww_0B/B$. Decompose $b \in
B$ into the product of an element of torus and of an unipotent
element: $b=tu$. Then $supp (Tx)=supp (Tbww_0B/B)=supp
(Tuww_0B/B)$. We represent $u$ as the exponent of an element of
the Lie algebra:
$$ u=\exp(\sum \limits _{\alpha_i \in \bar{\Delta}^{+}\subseteq
\Delta^{+}}c_i e_i)= \id + \frac{\sum c_i e_i}{1!}+\frac{(\sum c_i
e_i)^2}{2!}+\ldots,$$

\noindent where $e_i$ corresponds to the positive root $\alpha_i$;
$c_i\neq 0$, and $\bar{\Delta}^{+}\subseteq \Delta^{+}$ is a
subset of positive roots. We assume that  $u$ is written in the
normal form, so that $\bar{\Delta}^{+}\subseteq \Delta^{+}\cap
w\Delta^{+}$ ~\cite [28.4] {hum}.

{\vspace {1ex}}

As $supp(ww_0B/B)=ww_0\chi$ and the action of $T$ doesn't change
the support, from the exponent decomposition for  $u$ we get
$supp(Tx) \subset ww_0\chi+\sum \limits_{\alpha_i \in
\bar{\Delta}^{+}}\Bbb Z_{+} \alpha_i$.

{\vspace {1ex}}

We prove now the theorem describing the linear span of the support
in terms of the subsystem of positive roots
$\widetilde{\Delta}^{+}$. We assume that this root system is
saturated.

{\vspace {1ex}}

\begin{theorem} Let $w$  be an element of the Weyl group
such that
 $ww_0\in W^{st}_{\chi}$, and $\widetilde{\Delta}\subseteq\Delta$
be a saturated subsystem of roots  with $0\in (ww_0\chi+\sum
\limits_{\alpha_i \in \widetilde{\Delta}^{+}\cap w\Delta^{+}} \Bbb
Q_{+} \alpha_i).$ Then the support $supp(Tuww_0B/B)$ of the
$T$-orbit, where $u=\exp(\sum \limits _{\alpha_i \in
\widetilde{\Delta}^{+}\cap w\Delta^{+}}c_i e_i)$, contains zero
for almost all $c_i$.
\end{theorem}

\begin{proof} We begin with the following lemma.

{\vspace {1ex}}

\begin{lemma} Under the conditions of the previous theorem,
let $\alpha_0$
 be the highest root of the root system $\widetilde{\Delta}$. A strictly dominant
 weight is contained in the set of weights of the orbit  $M
\in (-C)^0\cap supp(Tuw_0B/B)$. Then the weight $\alpha_0+M$ is
also one of the  weights of the $T$-orbit.
\end{lemma}
\begin{proof} Let  $v_M$ be the vector of weight $M$, which
belongs to the weight decomposition of the vector corresponding to
the point of the orbit $Tx$. We may describe the vector $v_M$ by
the following formula:
$$v_M=\sum \limits_{M=w_0\chi+\sum \limits_{\alpha_i \in
\widetilde{\Delta}^{+}} a_i\alpha_i}
\frac{c_{\alpha_0}^{a_0}\ldots c_{\alpha_l}^{a_l}}{(\sum \limits_i
a_i)!}Sym(\underbrace{e_{\alpha_0},\ldots,
e_{\alpha_0}}_{a_0},\ldots,\underbrace{ e_{\alpha_l},\ldots,
e_{\alpha_l}}_{a_l})v_{w_0\chi},$$

\noindent where $Sym()$ denotes the sum over all permutations of
the products of the elements in the parentheses.

Note, that  the vector $v_M$ is nonzero for almost all values of
the coefficients $\{c_{\alpha_i}\}$ iff
 we may find the set of $\{a_i\}$ for which we have
$Sym(\underbrace{e_{\alpha_0},\ldots,
e_{\alpha_0}}_{a_0},\ldots,\underbrace{ e_{\alpha_l},\ldots,
e_{\alpha_l}}_{a_l})v_{w_0\chi}\neq 0$ (we denote this vector by
$v^0_{M}$). So to prove that the vector $v_{M+\alpha_0}$ is
nonzero  (for almost all values of $\{c_{\alpha_i}\}$) it is
sufficient to prove that
$Sym(e_{\alpha_0},\underbrace{e_{\alpha_0},\ldots,
e_{\alpha_0}}_{a_0},\ldots,\underbrace{ e_{\alpha_l},\ldots,
e_{\alpha_l}}_{a_l})v_{w_0\chi}\neq 0$.

{\vspace {1ex}}
 Let us show that the vector
$Sym(e_{\alpha_0},\underbrace{e_{\alpha_0},\ldots,
e_{\alpha_0}}_{a_0},\ldots,\underbrace{ e_{\alpha_l},\ldots,
e_{\alpha_l}}_{a_l})v_{w_0\chi}$ is proportional to
$e_{\alpha_0}v^0_{M}$. The element $e_{\alpha_0}$ commutes with
all $e_{\alpha_i}$, where $\alpha_i \in \widetilde{\Delta}^{+}$
($[e_{\alpha_i},e_{\alpha_0}]\subseteq \mathfrak
g_{\alpha_i+\alpha_0}=0$, as the weight $\alpha_i+\alpha_0$ is not
the root of the Lie algebra  $\mathfrak g$). That means
$e_{\alpha_0}$  can be taken out from the sign of symmetrization.

The proof of the lemma will be finished if we show that
$e_{\alpha_0}v^0_{M}\neq 0$. Consider the representation of the
${\mathfrak sl_2}\cong\langle
e_{\alpha_0},f_{\alpha_0},h_{\alpha_0}\rangle$, generated by
$v^0_{M}$. Hence  $e_{\alpha_0}v^0_M\neq 0$, as $M$ is strictly
antidominant and $v^0_M$ couldn't be the highest weight vector.
\end{proof}

\begin{center}
  \epsfig{file=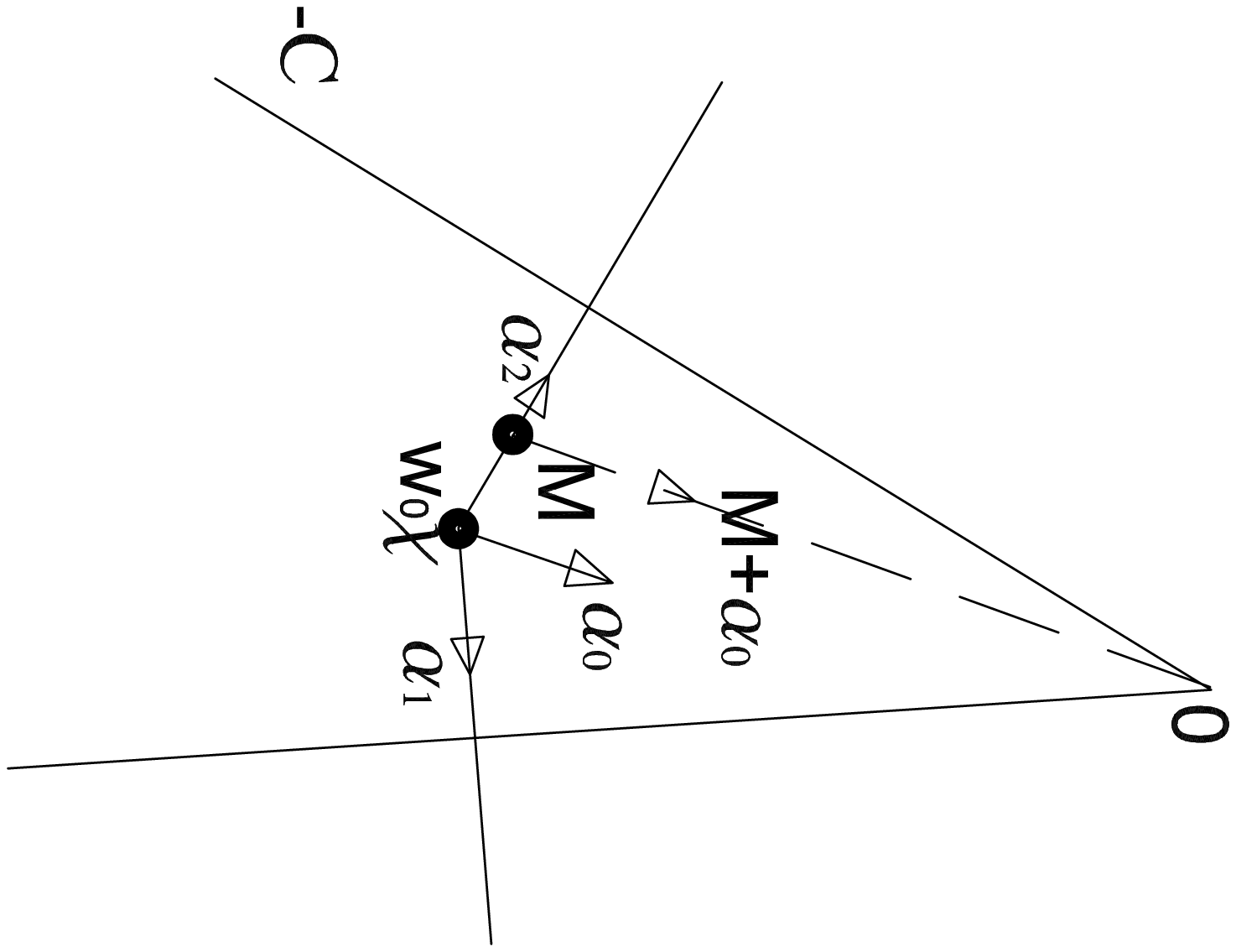,height=5cm,angle=90,clip=}
pic.2
\end{center}

\begin{remark} Instead of the claim that the weight $M$ is
antidominant we may require a weaker condition: $(M;\alpha_0)<0$.
Then we have $e_{\alpha_0}v^0_M\neq 0$. Indeed, from the condition
$(M;\alpha_0)<0$ it follows that $v^0_M$ couldn't be the highest
weight vector of the representation of $\mathfrak sl_2$ (the other
part of the proof remains the same).
\end{remark}

{\vspace {2ex}}

To prove next proposition  we  need to introduce some additional
notation. Denote by $M$ the rational weight and by
$\widetilde{\Delta}$ a subsystem of roots. Let
$H_M(\widetilde{\Delta}^+)$ be an affine cone $M+\sum
\limits_{\alpha_i \in \widetilde{\Delta}^+} \Bbb Q_+\alpha_i$ and
 $\delta H_M(\widetilde{\Delta}^+)$ its border. Sometimes we don't
 mention the root system in this notation that won't lead to the ambiguity.

{\vspace {2ex}}

\begin{proposition} Let $ww_0\chi$ is such weight that
$0=ww_0\chi+\sum \limits_{\alpha_i \in \Delta^+}c_i\alpha_i $,
where $c_i\geqslant 0$. Then every ray from $(-C)$ with the end in
zero point intersects the border  $\delta
H_{ww_0\chi}({\Delta}^+)$, and the intersection lies inside the
 $Pd_\chi$ polyhedron of the weights of the representation
$V(\chi)$.
\end{proposition}
\begin{proof} Consider the ray $l(t)=-t\sum b_i \pi_i$, where
$b_i\geqslant 0$ lying  $(-C)$ with the end in zero. We can
 describe $H_{ww_0\chi}$ by the inequalities
$(\pi_j,x)\geqslant(\pi_j,ww_0\chi)$. Denote by   $l(t_j)$ the
intersection point of the line containing the ray with the
hyperplane $j$, if it exists. (The line could be parallel to the
hyperplane. But it couldn't be parallel to all hyperplanes, as the
cone $H_{ww_0\chi}$ is solid. So we have at least one point of the
intersection.) For the proof of the first statement it is
sufficient to show that all $t_j\geqslant 0$. We may write the
equation of the intersection of the line and hyperplane $j$ in the
form $(\pi_j,l(t_j))-(\pi_j,ww_0\chi)=-t_j\sum \limits_{i} b_i
(\pi_i,\pi_j)+\sum \limits_{i}c_i(\alpha_i,\pi_j)=0$. So
$t_j=\frac {\sum \limits_{i}c_i(\alpha_i,\pi_j)}{\sum \limits_{i}
b_i (\pi_i,\pi_j)}\geqslant 0$ as $(\alpha_i,\pi_j)\geqslant 0$,
$(\pi_i,\pi_j)\geqslant 0$, $b_i,c_i\geqslant 0$. The minimal of
$t_j$ gives the intersection point of $l(t)$ and $\delta
H_{ww_0\chi}$.

\begin{center}
  \epsfig{file=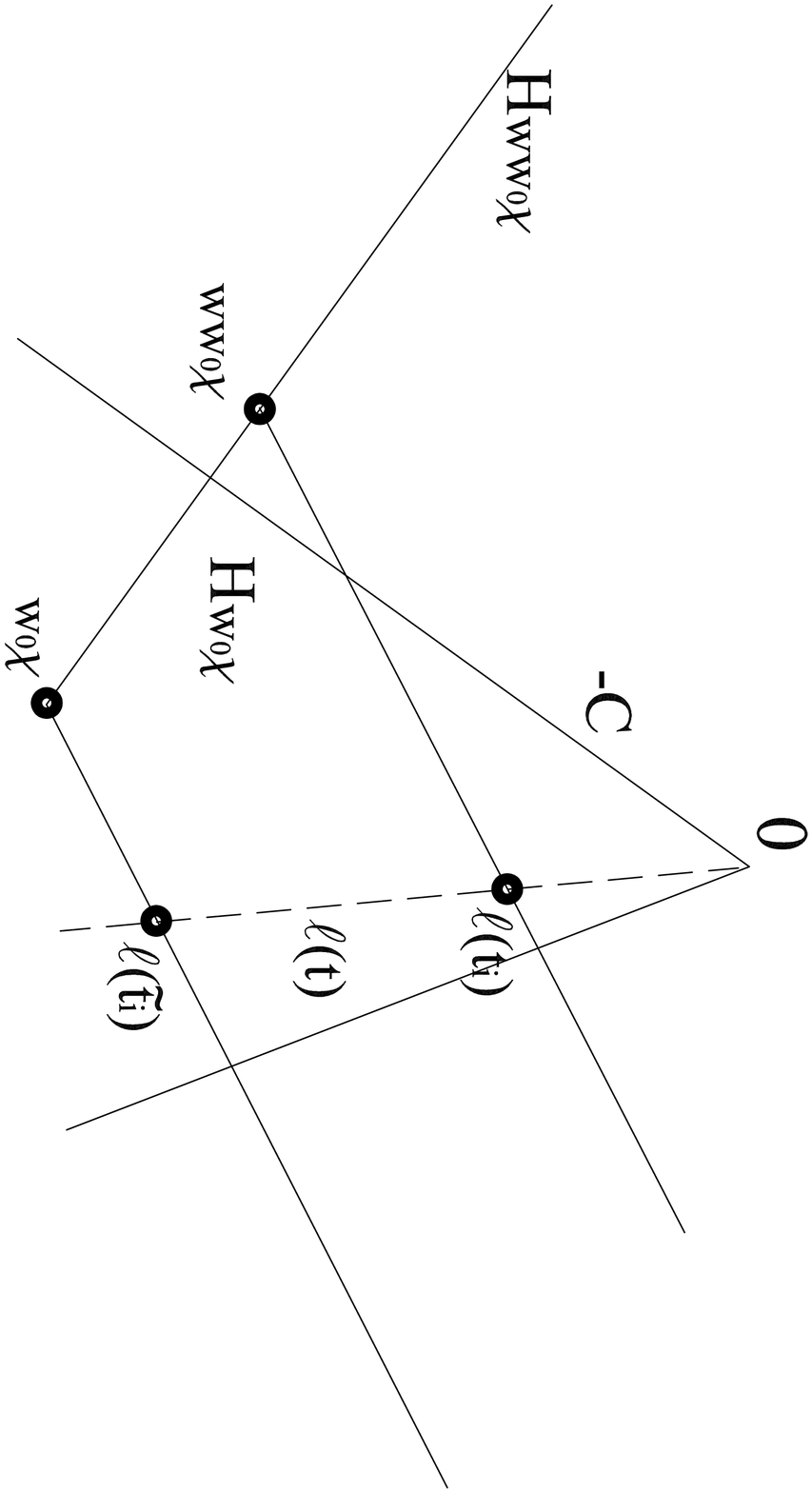,height=9cm,angle=90,clip=}
pic.3
\end{center}

We shall prove now that the intersection point lies inside the
weight polyhedron  $Pd_\chi$. It is sufficient to prove that the
intersection point belongs to the cone  $H_{w_0\chi}$ (as $P_\chi
\cap (-C)=H_{w_0\chi}\cap (-C) $). Consider the solutions $\tilde
t_j $ of the equations   $(\pi_j,l(\tilde t_j))=(\pi_j,w_0\chi)$ (
equations of the faces of the cone $H_{w_0\chi}$). As the faces of
the cone  $H_{w_0\chi}$ are parallel to the corresponding faces of
the cone $H_{ww_0\chi}$ it is sufficient to show that $\tilde t_j
\geqslant t_j$. We know that $ww_0\chi=w_0\chi+\sum \limits_{i}
d_i\alpha_i$, where $d_i\geqslant 0$. So we have $\tilde
t_j=\frac{\sum \limits_{i}d_i(\pi_j,\alpha_i)-\sum
\limits_{i}(\pi_j,ww_0\chi)}{\sum \limits_{i} b_i
(\pi_i,\pi_j)}\geqslant \frac{-\sum
\limits_{i}(\pi_j,ww_0\chi)}{\sum \limits_{i} b_i
(\pi_i,\pi_j)}=t_j$ as $(\pi_j,\alpha_i)\geqslant 0$ and
$d_i\geqslant 0$.
\end{proof} {\vspace {1ex}}

Now we prove the following lemma.

{\vspace {1ex}}

 \begin{lemma} Let  $ww_0\in
W^{st}_{\chi}$, $\alpha_0$ be the highest root of the saturated
root subsystem $\widetilde{\Delta}$,  for which we have
$ww_0\chi=-\sum \limits_{\alpha_i \in \widetilde{\Delta}^+}
c_i\alpha_i $, where $ c_i\geqslant 0$. Then $\alpha_0 \in
\widetilde{\Delta}^{+}\cap w\Delta^{+}$.
\end{lemma}

\begin{proof} The weight $w_0\chi$  is strictly antidominant.
Hence a root  $\alpha$ is positive  iff $\langle w_0\chi;\alpha
\rangle\leqslant 0$. As
 $ww_0\chi=-\sum
\limits_{\alpha_i \in \widetilde{\Delta}^+} c_i\alpha_i $, where
$c_i\geqslant 0$, we have $\langle ww_0\chi;\lambda
\rangle\leqslant 0$ for $\forall \lambda \in C_{\widetilde
\Delta}$. As the root $\alpha_0$ is the highest
$\widetilde{\Delta}$, it belongs to the closure of the Weyl
chamber $C_{\widetilde \Delta}$ and we may set $\lambda$ equal to
it. Consequently $\langle ww_0\chi;\alpha_0 \rangle=\langle
w_0\chi;w^{-1}\alpha_0 \rangle\leqslant 0$ that shows that the
root $w^{-1}\alpha_0$ is positive.
\end{proof}

{\vspace {2ex}}

We shall construct the piecewise-linear path from the zero point
to the point $ww_0\chi$. We proceed by the induction.

 {\vspace
{1ex}}

Let $\alpha_0$ be the highest root of the saturated subsystem
$\widetilde{\Delta}$. Besides we have  $ ww_0\chi=-\sum
\limits_{\alpha_i \in \widetilde{\Delta}^{+}\cap w\Delta^{+}} c_i
\alpha_i$, where $c_i\geqslant 0$.

Consider the ray with the end in zero point containing
$-\alpha_0$. As $\alpha_0$ is the highest root, this ray belongs
to the antidominant Weyl chamber $C_{\widetilde \Delta}$.
According to Proposition 2.4 it will intersect $\delta
H_{ww_0\chi}(\widetilde{\Delta}^{+})\cap
(-C_{\widetilde{\Delta}})$ in the  point denoted by $M_1$. So
$M_1$ belongs to the face of the cone
$H_{ww_0\chi}(\widetilde{\Delta}^{+})$. This face is the cone
 $H_{ww_0\chi}(\widetilde{\Delta}_1^{+})$
for a saturated root subsystem $\widetilde{\Delta}_1 \subset
\widetilde{\Delta}$. As  $M_1$ belongs to the antidominant Weyl
chamber it also belongs to the chamber
$-C_{\widetilde{\Delta}_1}$. So we may claim that
$M_1=-k_0\alpha_0$ for $k_0 \in \Bbb Q_+$.

Let us describe the induction step. On the i-th step we have the
sequence of the saturated root subsystems
${\Delta}^{+}\supseteq\widetilde{\Delta}_1^{+}\supseteq\ldots\supseteq
\widetilde{\Delta}_{i-1}^{+}\supseteq \widetilde{\Delta}_i^{+}$,
and also the sequence of the roots
$\{\alpha_0,\beta_1,\ldots,\beta_i\}$, where $\beta_j $, is the
highest root of the system  $\widetilde{\Delta}_j$. And we have
already constructed the sequence of weights
$\{0,M_1,M_2,\ldots,M_i\}$ such that  $M_j \in
H_{ww_0\chi}(\widetilde{\Delta}_j^{+})\cap
(-C_{\widetilde{\Delta}_j})$ and also
$M_j=M_{j-1}-k_{j-1}\alpha_{j-1}$, where $k_{j-1}\in \Bbb Q_+$.

Let us construct the weight $M_{i+1}$ and the root system
$\widetilde{\Delta}_{i+1}^{+}$. The weight  $M_i$ belongs to the
intersection $H_{ww_0\chi}(\widetilde{\Delta}_{i}^{+})\cap
(-C_{\widetilde{\Delta}_i})$. As the highest root  $\beta_i$ lies
in the Weyl chamber $C_{\widetilde{\Delta}_i}$, the intersection
point of the ray $M_i-t\beta_i$ (where $t\in\Bbb Q_+$) with the
end in $M_i$ and of the border
$H_{ww_0\chi}(\widetilde{\Delta}_{i}^{+})\cap
(-C_{\widetilde{\Delta}_i})$ belongs to  $\delta
H_{ww_0\chi}(\widetilde{\Delta}_{i}^{+})$. Denote by
$\widetilde{\Delta}_{i+1}^{+}\subset \widetilde{\Delta}_{i}^{+}$
the  saturated root subsystem corresponding to the face
$H_{ww_0\chi}(\widetilde{\Delta}_{i}^{+})$ which contains $M_i$.
Observe that  $M_{i+1} \in (-C_{\widetilde{\Delta}_{i+1}})$ and
that the number $k_{i+1}$ is rational as the weight $M_i$
 is rational and the cone $H_{ww_0\chi}(\widetilde{\Delta}_{i}^{+})$
is rational. Thus all the conditions are satisfied.

Our construction will be finished when $M_i=ww_0\chi$ for some
$i$. All weights $M_i$ are rational so one can find integer $N$
such that all $NM_i$ and also all $Nk_i$  are integer. If we take
the weight $Nww_0\chi$ instead  of $ww_0\chi$  (this will only
change the embedding $G/B \hookrightarrow \Bbb P(V(\chi))$ for
$G/B \hookrightarrow \Bbb P(V(N\chi))$, that don't affect the
statement of the theorem). The path constructed above will consist
of integer weights $NM_i$. Also these weights satisfy the
equalities $NM_{i+1}=NM_i-\tilde{k_i}\beta_i$, where $\tilde{k_i}$
 are integers. So we may suppose that  $M_i$ satisfy the
conditions  described above.

\begin{center}
  \epsfig{file=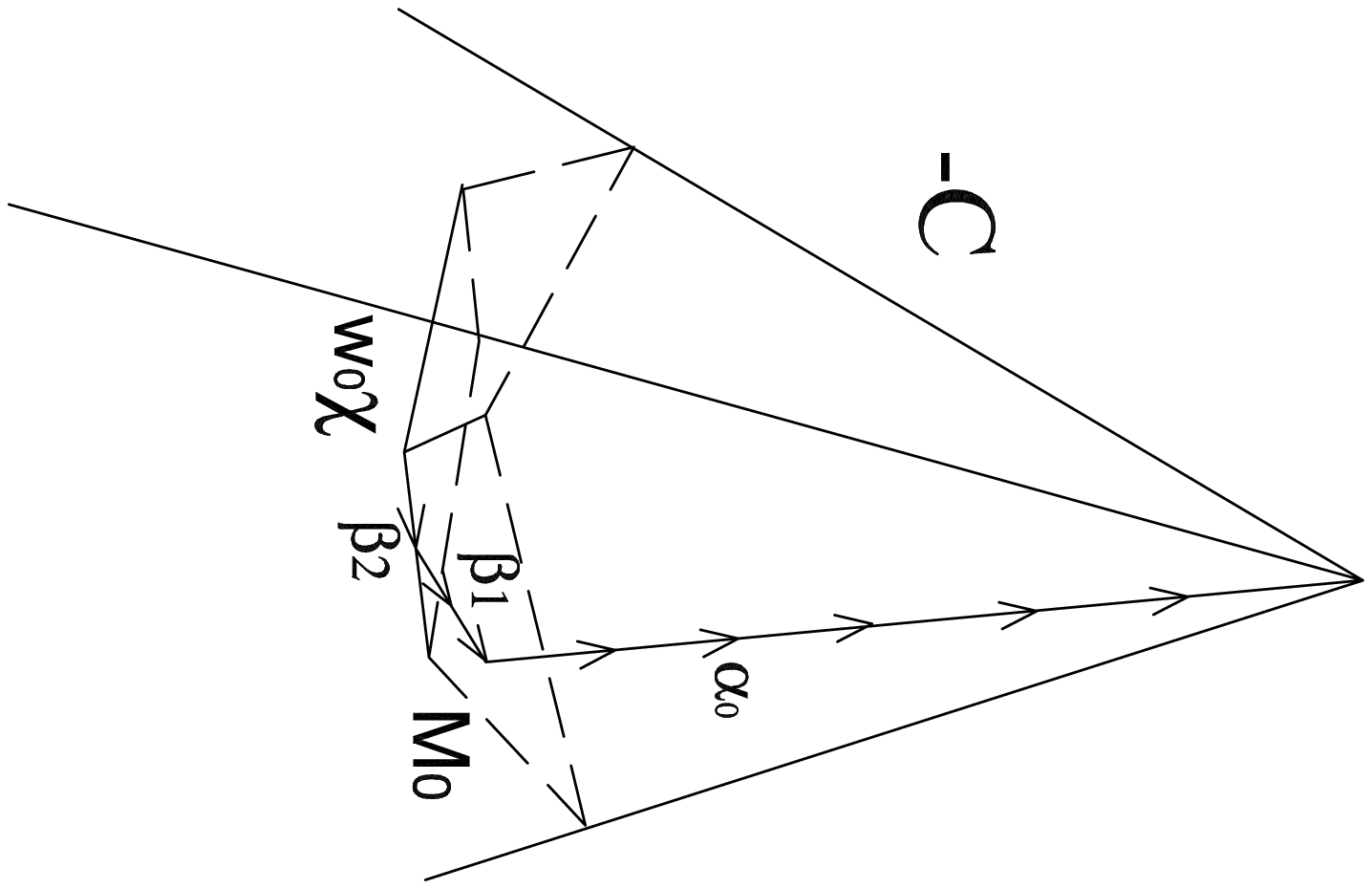,height=5cm,angle=90,clip=}
pic.4
\end{center}
{\vspace {1ex}}

We have to show now that all the weights $M_i$ belong to the
weights of the orbit $Tuww_0B/B$. In particular that will imply
that zero weight lies in the support of the orbit $Tuww_0B/B$. We
 argue by the induction. One may notice that $ww_0\chi$
belongs to the set of weights $Pd_\chi(Tuww_0B/B)$. Suppose that
the weight $M_{i+1}+l\beta_i \in Pd_\chi(Tuww_0B/B)$ for some $l$,
$0\leqslant l <k_i$ we must show that $M_{i+1}+(l+1)\beta_i \in
Pd_\chi(Tuww_0B/B)$.

In order to prove this we apply the sharper variant of  Lemma 2.2
(Remark 2.3) to the weight $M_{i+1}+l\beta_i$, root $\beta_i$ and
saturated root system $\widetilde{\Delta}_{i}$. But we have to
check that the conditions of the lemma are satisfied:
$(M_{i+1}+l\beta_i,\beta_i)< 0$ for every $0\leqslant l<k_i$
$(*)$, also we have to show that the root $\beta_i$ appears in the
exponent representing $u$ i.e. $\beta_i \in
\widetilde{\Delta}^{+}\cap w{\Delta}^{+}$ $(**)$. In the case
$\widetilde{\Delta}_{i}\varsubsetneq \Delta$ to apply the lemma to
the weight $M_{i+1}+l\beta_i$ and the system
$\widetilde{\Delta}_{i}$ it is required that  the formula for
$v_{M_{i+1}+l\beta_i}$ from  Lemma 2.2 contains only the terms
$e_{\alpha_j}$ with $\alpha_j \in \widetilde{\Delta}^{+}_i$ . But
this condition is satisfied because $M_{i+1}+l\beta_i \in
ww_0\chi+\sum \limits_{\alpha_j \in \widetilde{\Delta}_i^{+}}\Bbb
Q_{+} \alpha_j$ and the cone
$H_{ww_0\chi}(\widetilde{\Delta}_{i}^{+})$ is the face of
$H_{ww_0\chi}(\widetilde{\Delta})$ .

From the construction $M_{i+1}=-(k_i\beta_i+\ldots
+k_{1}\beta_{1}+k_0\alpha_0)$, where $\beta_j$, is the highest
root of $\widetilde{\Delta}_j^{+}$. Let $\beta_j$ be the highest
root of the system $\widetilde{\Delta}_j^{+}\supseteq
\widetilde{\Delta}_{j+1}^{+}\supseteq\ldots \supseteq
\widetilde{\Delta}_i^{+}$. Then $(\beta_j,\gamma)\geqslant 0$ for
every $\gamma \in  \widetilde{\Delta}_j^{+}$, in particular we
have $(\beta_j,\beta_m)\geqslant 0$ for every $m\geqslant j$. That
implies that $(\beta_j,\beta_i)\geqslant 0$ for every $j$, and
also $(\alpha_0,\beta_i)\geqslant 0$, as $\beta_i \in
\widetilde{\Delta}^{+}$ and $\alpha_0$ is the highest root of the
system $\widetilde{\Delta}^{+}$.

It follows from the above that
$(M_{i+1}+l\beta_i,\beta_i)=-((k_i-l)(\beta_i,\beta_i)+\ldots+k_{1}(\beta_{1},\beta_i)+k_0(\alpha_0,\beta_i))<
0$ that proves $(*)$.

Let us check $(**)$. By the construction
$-ww_0\chi=(k_l\beta_l+\ldots+k_{1}\beta_{1}+k_0\alpha_0)$ for
some $l$. One may notice that
$(-w_0\chi,w^{-1}\beta_i)=(-ww_0\chi,\beta_i)=k_i(\beta_l,\beta_i)+\ldots+k_{1}(\beta_{1},\beta_i)+k_0(\alpha_0,\beta_i)\geqslant
0$ as $(\beta_j,\beta_i)\geqslant 0$ and
$(\alpha_0,\beta_i)\geqslant 0$. From the inequality
$(-w_0\chi,w^{-1}\beta_i)\geqslant 0$ and the fact that the weight
$-w_0\chi$ is strictly dominant it follows that the root
$w^{-1}\beta_i$ is positive q.e.d.
\end{proof}

\begin{center}
  \epsfig{file=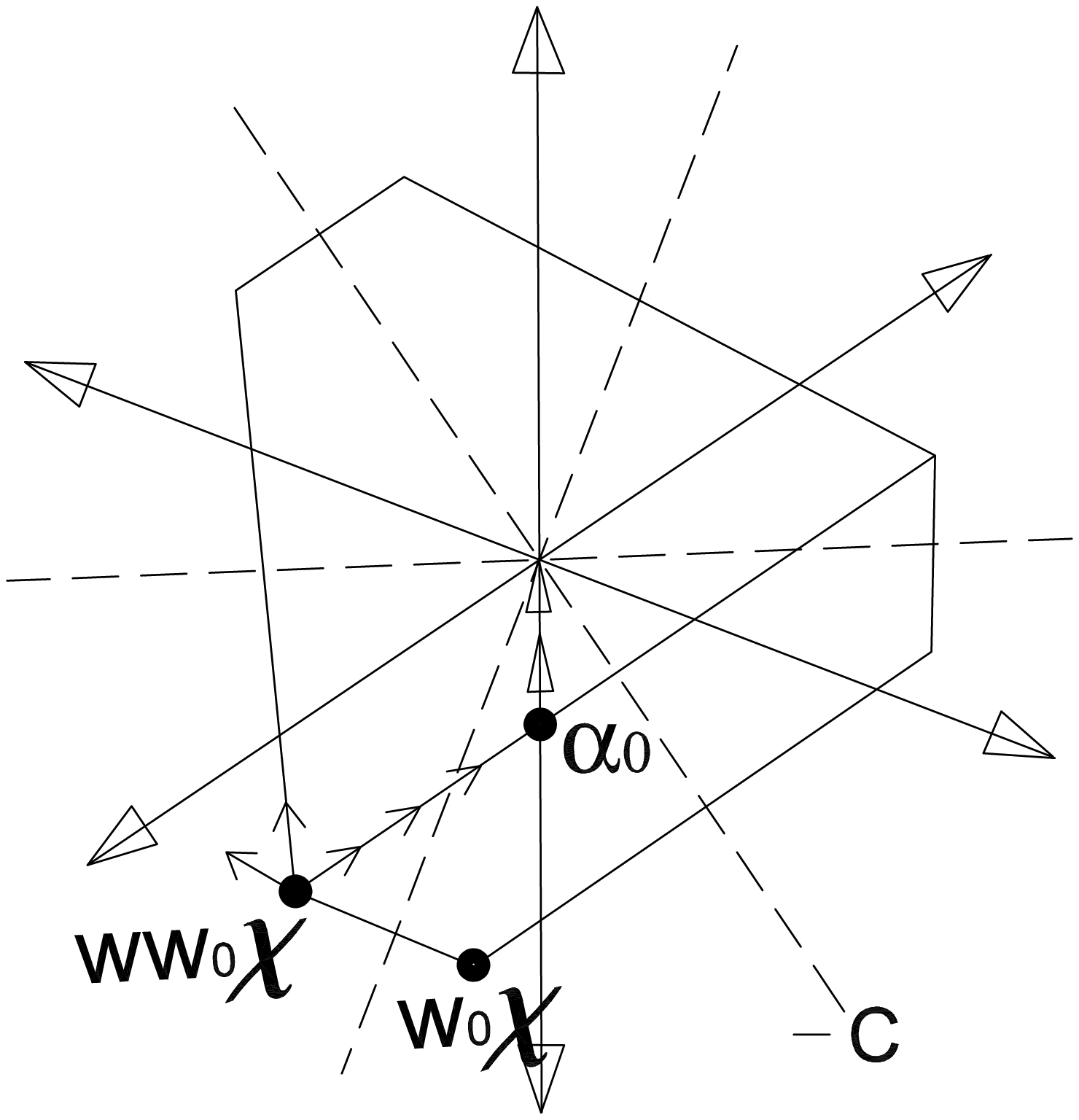,height=5cm,angle=90,clip=}
pic.5
\end{center}
{\vspace {1ex}}

{\vspace {2ex}}

\begin{corollary} Let $\chi$ be a strictly dominant weight
 and $w \in W$ be an element of Weyl group.  Then the following conditions are
 equivalent:

 (i) $0 \in ww_0\chi+\sum\limits_{\alpha_i \in {\Delta}^+}\Bbb
Q_+\alpha_i$,

 (ii) $0 \in ww_0\chi+\sum\limits_{\alpha_i \in {\Delta}^+\cap
w{\Delta}^+}\Bbb Q_+\alpha_i$.

These conditions imply that the support of the Schubert cell
$Bww_0B/B$ contains zero.
\end{corollary}

 \begin{proof}
 In the proof of the preceding theorem we constructed the piecewise-linear path from $ww_0\chi$
 to 0. The roots that appear in that path satisfy the condition $\alpha_i \in {\Delta}^+\cap w{\Delta}^+$.
 To construct this path we used only the fact that
 $0 \in ww_0\chi+\sum\limits_{\alpha_i \in {\Delta}^+}\Bbb
Q_+\alpha_i$. So the existence of such path implies $0 \in
ww_0\chi+\sum\limits_{\alpha_i \in {\Delta}^+\cap w{\Delta}^+}\Bbb
Q_+\alpha_i$.

Applying the preceding theorem we have $0 \in
Supp(Tuww_0B/B)\subset Supp(Bww_0B/B)$ for almost all values of
$c_i$ (where $c_i$ are the coefficients of $e_{\alpha_i}$
corresponding to the roots $\alpha_i \in {\Delta}^+\cap
w{\Delta}^+$ in the exponential decomposition of $u$).
\end{proof}

{\vspace {2ex}}

By means of this proposition one can describe the linear spans of
the supports for the semistable orbits.

{\vspace {1ex}}

\begin{proposition} Let  $Tuww_0B/B$ be a semistable
orbit from the cell $Bww_0B/B$. We claim that we can find the
orbit $T\tilde uww_0B/B$ from the same cell such that for $\tilde
u$ we have exponential decomposition $\tilde u=\exp(\sum \limits
_{\alpha_i \in \widetilde{\Delta}^{+}\cap w\Delta^{+}}c_i e_i)$,
where the subsystem $\widetilde{\Delta}$ is saturated and its
linear span is the same as the linear span of the support for
$Tuww_0B/B$.
\end{proposition}
\begin{proof} Consider the exponential decomposition
$u=\exp(\sum \limits _{\alpha_i \in {\widetilde{\Delta}}^{+}\cap
w\Delta^{+}}c_i e_i)$ (now we do not suppose that $c_i \neq 0$),
taking into account that  $u$ is chosen in the normal form i.e. $u
\in U\cap wU$. Let's consider the cone with the vertex $ww_0\chi$,
spanned by the roots $\alpha_i \in {\widetilde{\Delta}}^{+}\cap
w\Delta^{+}$ for which $c_i \neq 0$. Let's $\beta_i$ be the edges
of this cone. We show that the vectors
$e_{\beta_i}v_{ww_0\chi}\neq 0$ appear in the weight decomposition
of the vector from $V(\chi)$ corresponding to the point of the
orbit $Tuww_0B/B \subset \Bbb P(V(\chi))$.

Let us apply the element  $w^{-1}$ to the orbit. As this element
lies in $N_G(T)$ we shall get the  $T$-orbit $Tu'w_0B/B$ (where
$u'=w^{-1}uw$) and it will belong to the open cell (as $w_0\chi\in
supp(w^{-1}Tuww_0B/B)$ by  Lemma 1.1). Besides $\gamma \in
{\widetilde{\Delta}}^{+}\cap w\Delta^{+} $ is the root
corresponding to $e_{\gamma}$ in the exponential decomposition of
 $u$ iff the root $w^{-1}\gamma \in
{\Delta}^{+}\cap w^{-1}\widetilde{\Delta}^{+}$ corresponds to
$e_{w^{-1}\gamma}$ in the exponential decomposition of $u'$. So it
is sufficient to prove that $e_{w^{-1}\beta_i}v_{w_0\chi}\neq 0$.

Indeed, the weight $w^{-1}\beta_i \in {\Delta}^{+}$, $w_0\chi$ is
strictly antidominant so $v_{w_0\chi}$ couldn't be the highest
weight vector of the representation of ${\mathfrak sl_2}$ triple
$\langle
e_{w^{-1}\beta_i},f_{w^{-1}\beta_i},h_{w^{-1}\beta_i}\rangle$
generated by $v_{w_0\chi}$. That gives the claim.

The roots $\beta_i$ are the edges of the cone consequently the
weights $ww_0\chi+\beta_i$ belong to the support of the orbit.
Only the vectors $e_{\beta_i}v_{ww_0\chi}\neq 0$ could correspond
to them in the weight decomposition (we get this weight
decomposition by opening the brackets in the exponential
decomposition). The support of the orbit contains zero (the orbit
is semistable) so the linear span of the support coincides with
 the linear span of the roots $\langle \beta_1\ldots\beta_k\rangle_{lin}$. Consider the subsystem of roots
  $\widetilde{\Delta}=\Delta \cap \langle
\beta_1\ldots\beta_k\rangle_{lin}$. Applying the preceding theorem
to $\widetilde{\Delta}$ we  get the orbit of special type with the
same linear span of the support as the initial orbit
has.\end{proof}

{\vspace {2ex}}

 Thus to describe the linear spans of the supports for all
semistable orbits from the cell $Bww_0B/B$ it is sufficient to
describe only saturated root systems  $\widetilde{\Delta}$ for
which we have $0\in (ww_0\chi+\sum \limits_{\alpha_i \in
\widetilde{\Delta}^{+}\cap w\Delta^{+}} \Bbb Q_{+} \alpha_i)$.

{\vspace {3ex}}

\section{  The calculation of  ${\rm Pic}(X^{ss}/\!\!/T)$.}

{\vspace {2ex}}

According to Corollary 1.13, ${\rm Pic}_T(X^{ss}_{L_{\chi}})={\rm
Pic}_T(G/B)\cong {\rm Pic}(G/B)\times \Xi(T)$ (the last equality
follows from the fact that every line bundle is $T$-linearized and
every two $T$-linearizations differ by the character).

It is known that ${\rm Pic}(G/B)\cong\Xi(B)=\Xi(T)$
~\cite{popov},~\cite{kn}. The isomorphism  is constructed by the
following way: to a line bundle we associate the character of
action of $B$ on the fiber over $B$-stable point $eB/B$. And vica
versa, to a character $\chi$ we may associate the homogenious
bundle $G*_{B}k_{\chi}$ ( $k_{\chi}$ is a linear space, where $B$
acts by the character $\chi$) obtained as a quotient of the
$G\times k_{\chi}$ by $B$: $b(g,t)=(gb^{-1},\chi(b)t)$. Under this
isomorphism the cone of the very ample line bundles corresponds to
the internal points of the Weyl chamber.

$T$-linearized bundles may be written in the form
$G*_Bk_{\chi_0}\otimes k_{\chi_1}$. As a line bundle it is
isomorphic to  $G*_Bk_{\chi_0}$, but the torus action is twisted
by the character $\chi_1$.

Let  $\pi: X^{ss}_{L_\chi} \longrightarrow X^{ss}_{L_\chi}/\!\!/T
$ be the quotient morphism. It is well known that $ {\rm
Pic}(X^{ss}_{L_\chi}/\!\!/T)$ injects into ${\rm
Pic}_TX^{ss}_{L_\chi}$ by the map $\pi^*$ ~\cite{kn}. In the next
theorem we formulate the conditions when the bundle  $M\in {\rm
Pic}_TX^{ss}_{L_\chi}$ belongs to the subgroup $\pi^* {\rm
Pic}(X^{ss}_{L_\chi}/\!\!/T)$.

{\vspace {2ex}}

\begin{theorem} Let $\chi$ be a strictly dominant weight,
to which we associate the embedding $G/B \hookrightarrow \Bbb
P(V(\chi))$. Let $\{\widetilde{\Delta}^w_j \}$ be all saturated
root subsystems in  $\Delta$
 satisfying the condition
$0\in  ww_0\chi+\sum \limits_{\alpha_i \in
(\widetilde{\Delta}_j^{w})^+\cap w\Delta^{+}} \Bbb Q_{+}
\alpha_i$.
 Then the element $\mu=(\mu_0;\mu_1) \in ({\rm
Pic}(G/B)\times \Xi(T))\otimes{\Bbb Q}$ belongs to  $\pi^*{\rm
Pic}(X_{L_\chi}^{ss}/\!\!/T)\otimes \Bbb Q$ iff
$$ ww_0\mu_0+\mu_1\in
\bigcap \limits_{j }\langle\widetilde{\Delta}_j^{w}\rangle,$$ for
all $ww_0 \in W^{st}_{\chi}$. The rank of the Picard group ${\rm
Pic}(X_{L_\chi}^{ss}/\!\!/T)$  is equal to the dimension of the
linear space spanned by the points $\mu$ satisfying the conditions
described above.
\end{theorem}

\begin{remark} For the open cell the corresponding condition
could be simplified. Let $\{\widetilde{\Delta}_j\}$ be all
saturated root subsystems such that $w_0\chi
\in\langle\widetilde{\Delta}_j\rangle $. Then the orbits from the
open cell impose the following condition $$w_0\mu_0+\mu_1 \in
\bigcap \limits_{j}\langle\widetilde{\Delta}_j\rangle,$$
\end{remark}
 {\vspace {1ex}}

\begin{proof} Consider the following exact sequence from the work
~\cite{kn}:

$$ 1\longrightarrow {\rm
Pic}(X^{ss}/\!\!/T)\stackrel{\pi^*} {\longrightarrow} {\rm
Pic}_TX^{ss} \stackrel {\delta}{\longrightarrow} \prod \limits_{Tx
\subset X^{ss}} \Xi(T_x)$$

The last term is the  product over all orbits from $X^{ss}$ of
groups of characters of the stabilizers of  points in these orbits
( in the product it is sufficient to consider only the closed
orbits of $X^{ss}$). The map $\delta$ is the following: let us
take a line bundle $L \in {\rm Pic}_TX^{ss}$ and restrict it to
the orbit $Tx$. Then the stabilizer of the point in  $T_x$ will
act linearly on the fiber over $x$ by a character. This character
will be the image of $L$ in the component $\Xi(T_x)$ (denote this
map by $Pr_{\Xi(T_x)}:{\rm Pic}_TX^{ss}\longrightarrow \Xi(T_x)$),
so we have $\delta=\prod \limits_{Tx \subset
X^{ss}}Pr_{\Xi(T_x)}$.

{\vspace {1ex}}

Let  $\mu=\mu_0+\mu_1 \in {\rm Pic}(G/B)\times \Xi(T)$ be the
character defining line bundle $L_{\mu}$ from ${\rm Pic}_TX^{ss}$.

{\vspace {1ex}}

 Suppose we know  $\Ker(Pr_{\Xi(T_x)})$. Then
$ {\rm Pic}(X^{ss}/\!\!/T) \cong\Ker \delta = \bigcap \limits_{Tx
\subset X^{ss}}\Ker(Pr_{\Xi(T_x)})$.

{\vspace {1ex}}

As we calculate  the rank of the Picard group, it is sufficient to
consider only the action of one-parameter subgroups  $\lambda
:k^\times \rightarrow T_x$. Let  $v_x=v_{\tau_1}+\ldots +
v_{\tau_l}$ be the weight decomposition of $v_x$ for the embedding
$x=tuwB/B \in G/B \hookrightarrow \Bbb P(V(\chi))$, $\tau_i \in
\Xi(T)$.

 Consider the action of a one-parameter subgroup
$\lambda$ on $v_x$: $\lambda(t) v_x=t^{\langle
\lambda;\tau_0\rangle}v_x=t^{\langle \lambda;
\tau_1\rangle}v_{\tau_1}+\ldots +t^{\langle \lambda;
\tau_l\rangle} v_{\tau_l}$. Then  $\lambda(k^\times) \subset T_x$
iff $\langle \lambda;\tau_0\rangle=\langle \lambda;
\tau_1\rangle=\ldots=\langle \lambda; \tau_l\rangle$, or
equivalently  $\langle \lambda;\chi_i-\chi_j\rangle=0$ for all $
\chi_i,\chi_j \in supp(Tx \hookrightarrow \Bbb P(V(\chi)) )$. As
$0\in supp(Tx \hookrightarrow \Bbb P(V(\chi)) )$ this system of
equations is equivalent to $\langle \lambda;\chi_i\rangle=0$ for
$\forall \chi_i \in supp(Tx)$.

As in the proof of Proposition 2.7 we have $ww_0\chi \in supp
(Tx)$ and also $ww_0\chi+\alpha_i \in supp (Tx)$, where $\alpha_i$
are the roots appearing in the exponent decomposition of $u$ with
the nonzero coefficient. We suppose that $u$  is written in the
normal form. Then we have $\langle \lambda;\alpha_i\rangle=0$,
that means that  $T^0_x$ stabilizes $u$.

Let us calculate the character for the action of $T^0_x$ on the
fiber over $x$ of the bundle $G*_Bk_{\mu_0}\otimes k_{\mu_1}$. As
$T^0_x$ stabilizes $u$,  the character for the action of
 $T^0_x$ on the fiber over  $uww_0B/B$ is equal to the character for the action of $T^0_x$
  on the fiber over $ww_0B/B$, i.e.
$ww_0\mu_0+\mu_1$.

From the above it follows that  the character for the action of
 $T^0_x$ on the fiber will be trivial iff $\langle\lambda,ww_0\mu_0+\mu_1 \rangle=0$
 for every one-parameter subgroup from the stabilizer. So the condition for the
 line bundle $L_{\mu}$ to be in  $\Ker(Pr_{\Xi(T_x)})\otimes \Bbb Q$
can be rewritten in the form $ww_0\mu_0+\mu_1\in \langle supp(Tx
\hookrightarrow \Bbb P(V(\chi)) )\rangle $.

{\vspace {1ex}}

 Now we apply Proposition 2.7. The linear spans of the supports for the
 orbits of  special type from $X^{ss}\cap Bww_0B/B$
 are in correspondence with the subsystems of  positive roots
$\{\widetilde{\Delta}^{w+}_j\}$, satisfying the condition: $0\in (
ww_0\chi+\sum \limits_{\alpha_i \in
(\widetilde{\Delta}^{w}_j)^+\cap w{\Delta}^{+}} \Bbb Q_{+}
\alpha_i)$. {\vspace {1ex}}

Rewrite the conditions $\mu \in \Ker(\delta)$ in terms of the
weight  $\chi$:

$$ ww_0\mu_0+\mu_1\in \bigcap \limits_{j}\langle\widetilde{\Delta}_j^{w}\rangle.$$

 \end{proof}

{\vspace {2ex}}

\begin{example} Consider the case of general position; i.e.,
when the character  $w_0\chi$ does not belong to  linear spans of
 root subsystems in $\Delta$, those dimension is less than the
rank of the system $\Delta$. In this case there are no conditions
both on $\mu_0$ and $\mu_1$, so we have ${\rm rk} ({\rm
Pic}(X^{ss}_{L_\chi}/\!\!/T)) ={\rm rk}( {\rm
Pic}_T(X_{L_\chi}^{ss}))=2{\rm rk}G$ (in the case free of the
components of  type  $A_n$).
\end{example}
\begin{example} Consider the root system $B_4$. The simple roots
and the fundamental weights dual to them are the following:
$\alpha_1=\varepsilon_1-\varepsilon_2,\alpha_2=\varepsilon_2-\varepsilon_3,
\alpha_3=\varepsilon_3-\varepsilon_4,\alpha_4=\varepsilon_4$, and
 $\pi_1=\varepsilon_1, \pi_2=\varepsilon_1+\varepsilon_2,
\pi_3=\varepsilon_1+\varepsilon_2+\varepsilon_3,
\pi_4=\frac{1}{2}(\varepsilon_1+\varepsilon_2+\varepsilon_3+\varepsilon_4)
$. Consider the strictly dominant weight
$\chi=(\varepsilon_4-\varepsilon_3)+10(2\varepsilon_1
+\varepsilon_2+\varepsilon_3)=10\pi_1+\pi_2+8\pi_3+2\pi_4$. One
may check that  the only subspaces that contain this element are
the linear spans of the following root systems:  the first is
$\Delta_1$ of type $A_3$, generated by the roots
$\varepsilon_2-\varepsilon_3,\varepsilon_3-\varepsilon_4,\varepsilon_1+\varepsilon_4$,
which are the simple roots of this system (it also contain the
roots
$\varepsilon_2-\varepsilon_4,\varepsilon_1+\varepsilon_3,\varepsilon_1+\varepsilon_2$
as the positive roots);  the second is  $\Delta_2$  of  type
$A_1\oplus A_2$,  generated by
$\varepsilon_1,\varepsilon_2+\varepsilon_4,\varepsilon_3-\varepsilon_4$,
which are the simple roots in it (it also contains
$\varepsilon_2+\varepsilon_3$ as the positive root). The
intersection is  two dimensional space generated by the root
$\varepsilon_3-\varepsilon_4$ and the dominant weight
$2\varepsilon_1+\varepsilon_2+\varepsilon_3=\pi_1+\pi_3$.

We show that  ${\rm rk} ({\rm Pic}(X^{ss}_{L_\chi}/\!\!/T))=2$.

Let us write down the conditions of  Theorem 3.1 for the part of
the $\chi$-semistable elements
$w_0,s_{\alpha_1}w_0,s_{\alpha_2}w_0,s_{\alpha_3}w_0,s_{\alpha_4}w_0$
of the Weyl group:
$$w_0\mu_0+\mu_1 \in
\langle\Delta_1\rangle\cap\langle\Delta_2\rangle $$
$$w_0\mu_0+s_{\alpha_1}\mu_1 \in
\langle\Delta_1\rangle\cap\langle\Delta_2\rangle$$
$$w_0\mu_0+s_{\alpha_2}\mu_1 \in
\langle\Delta_1\rangle\cap\langle\Delta_2\rangle$$
$$w_0\mu_0+s_{\alpha_3}\mu_1 \in
\langle\Delta_1\rangle\cap\langle\Delta_2\rangle$$
$$w_0\mu_0+s_{\alpha_4}\mu_1 \in
\langle\Delta_1\rangle\cap\langle\Delta_2\rangle$$

 To prove the last four expressions one has to check that $ ( s_{\alpha_i}\chi\in\sum \limits_{\beta_j
\in (s_{\alpha_i}\Delta_1)^+\cap s_{\alpha_i}\Delta^+} \Bbb Q_{+}
\beta_j)$ and the same for the system $\Delta_2$. But
$(s_{\alpha_i}\Delta_1)^+\cap
s_{\alpha}\Delta^+=s_{\alpha_i}\Delta_1\cap
s_{\alpha_i}\Delta^+\cap \Delta^+$, so if we apply  $s_{\alpha_i}$
to  two last expressions we get $ ( \chi\in\sum \limits_{\beta_j
\in \Delta_1^+\setminus\{\alpha_i\}} \Bbb Q_{+} \beta_j)$. Indeed
 $s_{\alpha}$ is a simple reflection so $\Delta^+ \cap
s_{\alpha}\Delta^+=\Delta^+ \setminus{\{\alpha\}}$, that implies
$s_{\alpha_i}(s_{\alpha_i}\Delta_1\cap s_{\alpha_i}\Delta^+\cap
\Delta^+)=(\Delta_1^+\setminus\{\alpha_i\})$

We note that in the case $i=1,4$ the condition is satisfied
automatically since the roots $\alpha_1,\alpha_4$ do not belong to
the systems $\Delta_1$ è $\Delta_2$. In the case $i=2$ the second
condition is satisfied since $\alpha_2\notin \Delta_2$, and for
the first system it is satisfied since
$\chi=10(\varepsilon_1+\varepsilon_2)+9(\varepsilon_1+\varepsilon_3)+(\varepsilon_1+\varepsilon_4)
$.
 In the case
$i=3$ these conditions are true for $\Delta_1$ as:
$\chi=(\varepsilon_1+\varepsilon_4)+10(\varepsilon_2-\varepsilon_3)+19(\varepsilon_1+\varepsilon_3)$,
 and for the system $\Delta_2$ as: $\chi=20\varepsilon_1+9(\varepsilon_2+\varepsilon_3)+(\varepsilon_2+\varepsilon_4)$

 Subtract from the first expression the same expression for the root
 $s_{\alpha_i}$ so we get
$(w_0\mu_0+\mu_1)-(w_0\mu_0+s_{\alpha_i}\mu_1)=\frac{2(\alpha_i,\mu_1)}{(\alpha_i,\alpha_i)}\alpha_i\in
\langle\Delta_1\rangle\cap\langle\Delta_2\rangle$ for $i=1,2,3,4$.
As $\alpha_i \notin
\langle\Delta_1\rangle\cap\langle\Delta_2\rangle $ for $i=1,2,4$
we have $(\alpha_i,\mu_1)=0$. That implies $\mu_1=t\pi_3$. The
similar expression  for the element  $s_{\alpha_3}$ doesn't impose
any conditions since $\alpha_3 \in
\langle\Delta_1\rangle\cap\langle\Delta_2\rangle $. Let us show
that $\mu_1=0$. Then we would have that all the conditions imposed
by the semistable orbits would be the consequences of the
conditions imposed by the semistable orbits of the open cell i.e.
$w_0\mu_0 \in \langle\Delta_1\rangle\cap\langle\Delta_2\rangle$.
 Indeed if  $\mu_1=0$,  the conditions for $w \in
W^{st}_{\chi}$ can be rewritten in the form $w_0\mu_0 \in \bigcap
\langle\Delta_i\rangle$, where $\chi \in \Delta_i$, but all these
conditions are the consequences of the conditions for the open
cell, since the intersection is minimal for that. So we have
${\rm rk}( {\rm Pic}(X^{ss}_{L_\chi}/\!\!/T))=2$.

 To prove the equality $\mu_1=0$ consider the element
$w=s_{\alpha_3}s_{\alpha_4}$. First we prove that $ww_0$ is
$\chi$-semistable, and that it imposes the condition
$ww_0\mu_0+\mu_1 \in \langle w\Delta_1\rangle\cap\langle
w\Delta_2\rangle$ on $\mu$.  It is necessary and sufficient that
we have $w\chi\in\sum \limits_{\alpha_j \in (w\Delta_1\cap
w\Delta^+\cap \Delta^+)} \Bbb Q_{+} \alpha_j$ and the similar
condition for the root system $\Delta_2$. Apply  $w^{-1}$ to this
expression. Taking into account that  $w$ translates the roots
$(\varepsilon_3+\varepsilon_4)$ and $\varepsilon_4$ into negative
ones we get $\chi\in\sum \limits_{\alpha_j \in (\Delta_1\cap
\Delta^+\cap s_{\alpha_4}s_{\alpha_3}\Delta^+)} \Bbb Q_{+}
\alpha_j=\sum \limits_{\alpha_j \in
({\Delta^+_1}\setminus{\{\alpha_4,\varepsilon_3+\varepsilon_4\})}}
\Bbb Q_{+} \alpha_j$ and the similar condition for the system
$\Delta_2$. But both conditions are evidently satisfied since
$\varepsilon_4$ and $\varepsilon_3+\varepsilon_4$ do not belong to
either of the systems $\Delta_i$. So the conditions imposed by the
semistable orbits in the cell $Bww_0B/B$ can be rewritten as
$w_0\mu_0+w^{-1}\mu_1 \in
\langle\Delta_1\rangle\cap\langle\Delta_2\rangle$. Taking into
account that $\mu_1=t\pi_3$ we get
$w_0\mu_0+w^{-1}\mu_1=w_0\mu_0+ts_{\alpha_4}(\pi_3-\alpha_3)=
w_0\mu_0+t\pi_3-t(\varepsilon_3+\varepsilon_4) \in
\langle\Delta_1\rangle\cap\langle\Delta_2\rangle$. Subtract
$w_0\mu_0+t\pi_3 \in
\langle\Delta_1\rangle\cap\langle\Delta_2\rangle$  from the
previous expression. So we get $t(\varepsilon_3+\varepsilon_4) \in
\langle\Delta_1\rangle\cap\langle\Delta_2\rangle$, that is
possible only when $t=0$. Consequently $\mu_1=0$.
\end{example}

{\vspace {3ex}}

\end{document}